\documentstyle[12pt,draft]{article} \textheight 21cm \textwidth 15cm
\begin{document}

\title{ \Large {A Note on Characteristic Classes}
 \footnotetext{
\mbox{}\quad  \ \ {\it Subject classification}. \ 53C05, 55R25,
 57R20. }}
\author{ Jianwei Zhou }
\date{\small Department of Mathematics, Suzhou University, Suzhou
215006, P. R. China}

 \maketitle
\begin{abstract} This paper studies the relationship between  the sections
 and the Chern or Pontrjagin classes of a vector
bundle by the theory of connection. Our results are natural generalizations of the Gauss-Bonnet Theorem.
\end{abstract}

\baselineskip 15pt
\parskip 4pt

\vskip 1cm

\centerline{\bf 1. Introduction}

\vskip 0.3cm

Let $\pi: \ E \to M$ be an oriented Riemannian vector bundle with a Riemannian connection $D$, $\Omega$  its
curvature matrix. Then Euler characteristic class $e(E)$ of the bundle $E$ can be represented by $\mbox {Pf} (
\frac {-1}{2\pi} \Omega)$, where Pf is the Pfaffian polynomial. As is well-known, the Pontryagin and Chern classes
can be obtained from Euler classes. The characteristic classes are very important in the study of topology and
differential geometry.

In this paper, we  study the relationship between  the sections and the Chern or Pontrjagin  classes of a vector
bundle by the theory of connection. The results are natural generalizations of the Gauss-Bonnet Theorem which
concerns the relationship among the Euler class of the tangent bundle, the tangent vector fields and the
Euler-Poincar\'{e} characteristic number of the manifold.

As is well-known,  the top Chern class of a complex vector bundle $E_{\bf C}$ and the Euler class of its
realization vector bundle $E_{\bf R}$ are the same. This can be proved by splitting principle, see [1], p.273,
[4], p.115 or [10], p.158. In \S 2, we give a direct proof of this fact. Then we state some known results about
the characteristic classes which are needed in \S 3.

In \S 3, we study the relationship between the sections of the vector bundle and the Chern or Pontrjagin  classes
of the bundle. Using the transgression, we show that the Chern and the Pontrjagin classes  can be represented by
cycles in homology of the base manifolds by Poincar\'{e} dual. These cycles are determined by the generic sections
of the vector bundles.

In the following, we assume that the base manifolds of vector bundles are all compact and oriented.

\vskip 1cm

\centerline{ \bf 2.  Preliminaries}

\vskip 0.3cm

The complex Euclidean space ${{\bf C}}^n$ is naturally isomorphic to a real Euclidean space ${\bf R}^{2n}$. The
isomorphism can be given by
$$\sum z_ie_i \mapsto \sum x_ie_i +\sum y_ie_{n+i}, \ \
z_i=x_i+\sqrt {-1} y_i, \ i=1,\cdots,n,$$ where $e_1,\cdots,e_n$ is an unitary basis of ${{\bf C}}^n$. The basis
$e_1,e_{n+ 1},\cdots, e_n,e_{2n}$ of ${\bf R}^{2n}$  also gives an orientation on ${\bf R}^{2n}$. For any matrix
$C=(C_{ij})\in \underline {so}(2n)$ the Lie algebra of $SO(2n)$, let
$$T=(e_1,e_{n+ 1},\cdots,e_n,e_{2n})\wedge C (e_1,e_{n+
1},\cdots,e_n,e_{2 n})^t.$$ The pfaffian $\mbox {Pf} (C)$  is defined by $$\mbox {Pf} (C)e_1\wedge e_{n+1}\wedge
\cdots\wedge e_n\wedge e_{2n} = \frac 1{2^nn!} T^n.$$

Let $U(n)$ be the unitary group and $\underline {u}(n)$  its Lie algebra, any element of $\underline {u}(n)$ can
be represented by $ A+\sqrt {-1}B$, where $A,B$ are real matrices. The canonical representation $U(n) \to SO(2n)$
induces a representation between their Lie algebras.  With the oriented bases $e_1,\cdots,e_n$ and
$e_1,\cdots,e_n,e_{n+1},\cdots,e_{2n}$ on ${{\bf C}}^n$ and ${\bf R}^{2n}$ respectively, the map $\underline
{u}(n) \to \underline {so}(2n)$ can be represented by $$ A+\sqrt{-1} B \mapsto \left(
\begin{array}{rr} A & B \\ -B & A \end{array} \right), \ \ A^t=-A,\ B^t=B.$$
Denote $C$ the matrix obtained by rearrange the rows and columns of $\left(
\begin{array}{rr} A & B \\ -B & A \end{array} \right)$
 according to the oriented basis $e_1,e_{n+
1},\cdots,e_n,e_{2n}$ of ${\bf R}^{2n}$.

\noindent {\bf Lemma 2.1.} \ $\mbox {Pf} (C) = \det (-\sqrt{-1}(A+\sqrt {-1}B)).$

\noindent {\it Proof}.  \ It is easy to see that $T$ can also be represented by $$T =(e_1,\cdots,e_n,e_{n+
1},\cdots,e_{2n})\wedge \left(
\begin{array}{rr} A & B
\\ -B & A
\end{array} \right)(e_1,\cdots,e_n,e_{n+
1},\cdots,e_{2n})^t.$$ Let $g_i=e_i-\sqrt{-1}e_{n+i}, \ g_{n+ i}=e_i+\sqrt{-1}e_{n+i}$ and $X=\frac {\sqrt 2}2
\left(
\begin{array}{cc} I & \sqrt {-1}I \\ \sqrt {-1}I &  I
\end{array} \right)\in U(2n)$. By
$$(e_1,\cdots,e_n,e_{n+
1},\cdots,e_{2n})X=\frac {\sqrt 2}2(g_{n+1},\cdots,g_{2n},\sqrt{-1}g_{1},\cdots,\sqrt{-1}g_{n}),$$
$$X^{-1}(e_1,\cdots,e_n,e_{n+
1},\cdots,e_{2n})^t=\frac {\sqrt 2}2(g_{1},\cdots,g_{n},-\sqrt{-1}g_{n+1},\cdots,-\sqrt{-1}g_{2n})^t,$$
$$X^{-1}\left( \begin{array}{rr} A & B
\\ -B & A
\end{array} \right)X=\left( \begin{array}{rr} A+\sqrt{-1}B & {}
\\ {} & A-\sqrt{-1}B
\end{array} \right),$$
\begin{eqnarray*} &&
(g_{n+1},\cdots,g_{2n})\wedge(A+\sqrt{-1}B)(g_{1},\cdots,g_{n})^t\\
& & = (g_{1},\cdots,g_{n})\wedge(A-\sqrt{-1}B)(g_{n+1},\cdots,g_{2n})^t,\end{eqnarray*} we have $ T =
(g_{n+1},\cdots,g_{2 n})\wedge(A+\sqrt{-1}B)(g_1,\cdots,g_n)^t.$ Then
\begin{eqnarray*} &  T^n & =n! \det (
A+\sqrt {-1}B ) g_{n+1}\wedge g_1\wedge\cdots \wedge g_{2 n}\wedge
g_n \\
& & =n!\left(-2\sqrt {-1} \right)^n\det ( A+\sqrt {-1}B ) e_1\wedge e_{n+1}\wedge \cdots\wedge e_n\wedge e_{2n}.
\end{eqnarray*}
Hence $\mbox {Pf} (C)= \left(-\sqrt {-1} \right)^n\det ( A+\sqrt {-1}B ). \ \ \ \Box $

Let $\pi: \ E_{{\bf C}} \to M$ be a Hermitian vector bundle with the fibre ${{\bf C}}^n$, \ $D_{{\bf C}}$  a
Hermitian connection on $E_{{\bf C}}$. The bundle $E_{{\bf C}}$ naturally determines a real Riemannian vector
bundle  $\tau: \ E_{\bf R} \to M$ with fibre ${\bf R}^{2n}$ and  a Riemannian connection $D_{\bf R}$. If
$s_1,\cdots,s_n$ is an unitary basis for the sections of $E_{{\bf C}}$ over a open set $U\subset M$, then
$s_1,s_{n+1} =\sqrt {-1}s_1,\cdots,s_n,s_{2n} =\sqrt {-1}s_n$ form an orthonormal basis for the sections of
$E_{\bf R}$ over $U$, see [10], p.155. If $D^2_{\bf C}s_i=\sum\limits_{j=1}^n \widetilde \Omega_{ij} s_j, \ \
\widetilde \Omega_{ij} = \Omega_{ij}+ \sqrt {-1} \Omega_{i,n+j}$ the curvature forms of connection $D_{{\bf C}}$,
we have
$$D^2_{\bf R} s_i = \sum \Omega_{ij}s_j +\sum \Omega_{i,n+ j}s_{n+
j},$$
 $$D^2_{\bf R}s_{n+i} = - \sum \Omega_{i,n+ j}s_{j}+\sum \Omega_{ij}s_{n+ j}. $$

Denote $\Omega_{E_{\bf C}}$ and $\Omega_{E_{\bf R}}$ the curvature matrices of $D_{\bf C}, D_{\bf R}$
respectively. By Lemma 2.1, we have

\noindent {\bf Corollary 2.2.} \ {\it The top Chern class of bundle $E_{\bf C}$ and Euler class of $E_{\bf R}$
represented by $\Omega_{E_{\bf C}}$ and $\Omega_{E_{\bf R}}$ respectively are the same, that is, $\det (\frac
{\sqrt {-1}}{2\pi} \Omega_{E_{{\bf C}}})= \mbox {Pf} (\frac {-1}{2\pi}\Omega_{E_{\bf R}}).$ }

From the vector bundle $\pi: \ E_{{\bf C}} \to M$, we can construct fibre bundles $\pi_i: \ V(E_{{\bf C}},i)\to M,
\ i=1,\cdots,n$. For any $p\in M$, the fibre $\pi_i^{-1}(p)$  is a complex Stiefel manifold formed by all unitary
$i$-frames  on $\pi^{-1}(p)$. For each $i$, we have an induced bundle $\pi_i^*E_{{\bf C}} \to V(E_{{\bf C}},i)$
which can be decomposed by $\pi_i^*E_{{\bf C}} = {\cal E}_i \oplus F_{n-i}$. The fibre of $F_{n-i}$ over
$(s_1,\cdots,s_i)\in V(E_{{\bf C}},i)$ is the orthogonal complement of $s_1,\cdots,s_i$ in vector space
$\pi^{-1}(p)$, the bundle ${\cal E}_i$ is  trivial. Then we have the following commutative diagram
$$\begin{array}{cccccccclc} {} & {} & F_1   & \longrightarrow &
\cdots & \longrightarrow & F_{n-1} & \longrightarrow & F_n=E_{{\bf C}} \\
{} & {} & \downarrow & {}  & {} &{}    & \downarrow
& {} & {} \downarrow \\
V(E_{{\bf C}},n) & \stackrel {\alpha_n} \longrightarrow & V(E_{{\bf C}},n-1) &  \stackrel{\alpha_{n-1}}
\longrightarrow  & \cdots & \longrightarrow & V(E_{{\bf C}},1) & \stackrel{\alpha_1}
 \longrightarrow & M. \end{array} $$
The maps in the diagram are all defined naturally. By the theory of characteristic class and Corollary 2.2,
$$\pi_{n-i}^*c_i(E_{\bf C})= c_i(\pi_{n-i}^*E_{{\bf C}})
=c_i({\cal E}_{n-i} \oplus F_i)=c_i(F_i)= e({F_i}_{\bf R}).$$ The map $\alpha_j: \ V(E_{{\bf C}},j) \to V(E_{{\bf
C}},j-1)$ defines a fibre bundle with the fibre $S^{2n-2j+1}$. As in [10], \S 14, applying the   Gysin sequence to
the vector bundle $F_{n-j+1}\to V(E_{{\bf C}},j-1)$, we know that the pullback map
$$\alpha_j^*: \ H^k(V(E_{{\bf C}},j-1),Z) \to H^k(V(E_{{\bf C}},j),Z)$$ is an
isomorphism for any $ k< 2n-2j+1$. Since $\pi_{n-i} =\alpha_1\cdots \alpha_{n-i}: \ V(E_{{\bf C}},n-i)\to M$, the
maps
$$\pi_{n-i}^*: \ H^k(M,Z) \to H^k(V(E_{\bf C},n-i),Z), \ \ \
k<2i+1,$$ are all isomorphism.

\noindent {\bf Proposition 2.3.} \ $c_i(E_{{\bf C}})=\pi_{n-i}^{*-1} c_i(F_i)=\pi_{n-i}^{*-1}e({F_i}_{\bf R}), \
i=1,\cdots, n$.

Then $c_i(E_{{\bf C}})=0$ if and only if the bundle $F_i \to V(E_{{\bf C}},n-i)$ has a non-zero section.

For the  real vector bundle $\pi: \ E\to M$, we  can also construct fibre bundles $\pi_i: \ V(E,i)\to M, \
i=1,\cdots,n=$ rank $E$. For any $p\in M$, the fibre $\pi_i^{-1}(p)$  is a Stiefel manifold formed by all
orthonormal frames on $\pi^{-1}(p)$. For any $i$, we have a pullback vector bundle $\pi_i^*E = {\cal E}\oplus
\widetilde F_{n-i}\to V(E,i)$, where ${\cal E}$ is a trivial bundle of rank $i$.

\noindent {\bf Proposition 2.4.} \ {\it Assume the vector bundle $E$ is oriented, then the vector bundles
$\widetilde F_{n-k}\to V(E,k)$ are all oriented. We have $p_i(E) =\pi_{n-2i}^{*-1}[e(\widetilde F_{2i})\cdot
e(\widetilde F_{2i})],$ where $p_i(E)$ is $i$-th Pontrjagin class. }

\noindent {\it Proof}. \ Similar to Proposition 2.3, we have $p_i(E) =\pi_{n-2i}^{*-1}p_i(\widetilde F_{2i}).$ By
[10], Corollary 15.8, $p_i(\widetilde F_{2i})=e(\widetilde F_{2i})\cdot e(\widetilde F_{2i})$. \ \ \ $\Box$

\newpage

\vskip 1cm

\centerline{\bf  3. The Transgression and the Poincar\'{e} dual} \vskip 0.3cm

In [2],[3], Chern gave an elegant proof of the Gauss-Bonnet theorem and introduced the concept of the
transgression for the characteristic classes. Let $\pi:  \ E\to M$ be an oriented Riemannian vector bundle with
rank $2n,  \ p: \ S(E)=V(E,1) \to M$ be the associated sphere bundle. The induced bundle $p^*E \to S(E)$ can be
decomposed as $p^*E =\widetilde F_{2n-1}\oplus{\cal E}.$ Then
$$p^*e(E)=e(\widetilde F_{2n-1}\oplus{\cal E})=0 \ \ \mbox{in } \
H^{2n}(S(E),Z).$$ Let $D$ be a Riemannian connection on $E$ and $p^*D$ the pull back connection on $p^*E$. Let
$e_1,\cdots,e_{2n-1},e_{2n}$ be oriented orthonormal frame fields on $p^*E, \ {\cal E}$ be generated by $e_{2n}\in
S(E)$. Define another connection $\widetilde D$ on $p^*E$:
$$\widetilde D e_\alpha=\sum
\widetilde \omega_\alpha^\beta e_\beta, \ \alpha,\beta=1,\cdots,2n-1, \ \ \widetilde D e_{2n}=0,$$ where
$\widetilde \omega_\alpha^\beta= p^*\omega_\alpha^\beta$ are defined by $D e_\alpha=\sum  \omega_\alpha^\beta
e_\beta +\omega_\alpha^{2n}e_{2n}.$ Let $p^*\Omega$ and $\widetilde \Omega$ be the curvature matrices of the
Riemannian connections $p^*D$ and $\widetilde D$ on $p^*E$ respectively. Then $e(E)$ can be represented by
$e(\Omega)= \mbox {Pf} (\frac {-1}{2\pi}\Omega)$ and $e(\widetilde \Omega)=0$ on $p^*E$. By Chern-Weil methods,
there is a $2n-1$ form $\eta$ on $S(E)$ such that
$$p^*e(\Omega) = - d\eta, \ \
\eta =\frac {1}{(-2\pi)^n}\int_0^1 \mbox {Pf}(\widetilde \omega-p^*\omega,\Omega_t,\cdots,\Omega_t)dt,$$ where
$\Omega_t$ the curvature matrix of the connection $p^*D+t(\widetilde D-p^*D)$. Restricting $\eta$ to each fibre of
$S(E)\to M$ is the volume form of the fibre. For the computation of $\mbox {Pf}(\widetilde
\omega-p^*\omega,\Omega_t,\cdots,\Omega_t)$, see [7], p.297.

When $E=TM$ is the tangent bundle of Riemannian manifold $M$, the form $\eta$ is the same as Chern obtained in
[3].

Let $\rho: \ [0,+\infty) \to R$ be a smooth function, $\rho (r)=-1$ for $r\in [0,1]$, \ $\rho(r)=0$ for $r\geq 2$.
Then the $2n$-form $\Phi= d(\rho(|e|)\tau^*\eta)$ is a Thom form on $E$, where $|e|$ is the norm of $e\in E$ and
$\tau: \ E- M\to S(E)$ is the projection, $e\in E-M,  \ \tau(e)=e/|e|$. For proof, see [1], p. 132, Proposition
12.3. For the construction of Thom form, see also Mathai and Quillen [9].

Similarly, for the complex vector bundle $\pi: \ E_{{\bf C}}\to M$ defined in \S 2, we have

\noindent {\bf Theorem 3.1.} \ {\it For any $i=1,\cdots,n=$ rank $E_{{\bf C}}$, there is  a $2i-1$ form $\eta_i$
on $V(E_{\bf C},n-i+1)$ such that $\pi_{n-i+1}^*c_i(\Omega_{E_{{\bf C}}})=-d\eta_i.$ }

The theorem follows from Proposition 2.3 and the result on the Euler classes. Using the transgression form
$\eta_i$, we can construct a Thom form $\Phi_i$ for vector bundle $F_i$.

Let $s_1,\cdots,s_{n-i+1}$ be sections of the Hermite bundle $E_{{\bf C}}$ which are linearly independent on $M-
Z$. Assume that $Z$ is a set of submanifolds of $M$. From these sections, we have a section $\tilde s: \ M-Z \to
V(E_{{\bf C}}, n-i+1)$. Then on the subset $M-Z$, we have
$$c_i(\Omega_{E_{{\bf C}}})=-d(\tilde s^*\eta_i).$$ Let $U_\varepsilon$ be a
$\varepsilon$-neighborhood of $Z$ in $M$. For any closed $m-2i$ form $\xi$ on $M, \ m=\dim M$, we have $$\int_M
c_i(\Omega_{E_{\bf C}})\wedge \xi =\lim\limits_{\varepsilon \to 0} \int_{\partial U_\varepsilon} \tilde
s^*\eta_i\wedge\xi.$$ The left-hand side of this equation is independent of the choice of the sections of the
bundle $E_{{\bf C}}$. This equation is useful  for our understand the relationship between the characteristic
classes and the sections of the vector bundles as we know for the Euler classes.

\noindent {\bf Theorem 3.2.} \  {\it Let $s_1,\cdots,s_{n-i_1+1}$ be sections on $E_{{\bf C}}$ which are linearly
independent on $M-Z$, where $Z$ is a subset  of $M$. If there is a nonzero Chern number $a=\int_M
c_{i_1}(\Omega_{E_{{\bf C}}})\cdots c_{i_k}(\Omega_{E_{{\bf C}}}), \ i_2\geq\cdots \geq i_k, \ k > 1$. Then the
set $Z$ cannot be contained in any submanifold of $M$ with dimension less  than $2i_2$. }

Similar result holds on the real vector bundles.

\noindent {\it Proof}. \ If $Z$ is contained in a submanifold $N$ of dimension less then $2i_2$, let $U_1\subset
U_2$ be open neighborhoods of $N$ such that $N$ is a deformation retract of $U_2$. Then  on $U_2$, we have
$E_{{\bf C}}|_{U_2} = E_1\oplus E_2$ and $E_1$ is trivial with rank $
> n-i_2$. Then we can construct a connection $D_{{\bf C}}$ on $E_{{\bf C}}$ such
that $c_{i_2}(\Omega_{{\bf C}})|_{U_1}=0$, where $\Omega_{{\bf C}}$ is the curvature matrix of $D_{{\bf C}}$.
Hence
$$\int_M
c_{i_1}(\Omega_{E_{{\bf C}}})\cdots c_{i_k}(\Omega_{E_{{\bf C}}}) =\lim\limits_{\varepsilon \to 0} \int_{\partial
U_\varepsilon} \tilde s^*\eta_{i_1}\wedge c_{i_2}(\Omega_{E_{{\bf C}}})\cdots c_{i_k}(\Omega_{E_{{\bf C}}})=0,$$
 contradict to the fact of $a
\not= 0$. \ \ \ \ $\Box$

For example, let ${\bf C}P^n$ be the complex projective space and $T_c{\bf C}P^n$ its holomorphic tangent space.
It was proved in [10] that all Chern numbers of $T_c{\bf C}P^n$ are nonzero. Let $i_1=n-k, i_2=k$ in this case.
Then for any submanifold $N$ of ${\bf C}P^n$, $\dim_{\bf R} N <2k$, there do not exist vector fields
$X_1,\cdots,X_{k+1}\in \Gamma (T_c{\bf C}P^n)$ such that they are linearly independent on ${\bf C}P^n-N$.

In the following we give a geometric proof of the fact that the Chern   classes can be represented by some
submanifolds of the base manifold of the bundle which is the Poincar\'{e} dual.

As in [6], Chapter 3, let $\sigma =\{s_1,\cdots,s_{n}\}$ be a set of $C^\infty$ sections of bundle $E_{{\bf C}}$
and the degeneracy set $D_i(\sigma)$ be defined by
$$D_i = D_i(\sigma) =\{ p\in M \ | \ s_1(p)\wedge\cdots\wedge s_i(p) =0\},
 \ \ i=1,\cdots,n.$$
We say that $\sigma$ is generic if, for each $i, \ s_{i+1}$ intersects the subspace of $E_{{\bf C}}$ spanned by
$s_1,\cdots,s_{i}$ transversely, so that, $D_{i+1}$ is, away from $D_i$, a submanifold of dimension $m-2n+2i$.
Thus sections $s_1,\cdots,s_{i+1}$ are linearly independent everywhere if $m+2i< 2n$. We can give each $N_{i+1}=
D_{i+1}-D_i$ an orientation defined naturally. Then $D_{i+1}$ represents a cycle in homology, called the
degeneracy cycle of the sections $\sigma$. In a neighborhood of a point $p\in N_{i+1}$, complete the sections
$e_1=s_1,\cdots,e_i=s_i$ to a frame for $E_{{\bf C}}$, and write $s_{i+1}=\sum\limits_{j} f_je_j, \ f_j=
f_{j,1}+\sqrt {-1}f_{j,2}$, where $f_{j,1},f_{j,2}$ are real functions. $N_{i+1}$ is then locally given by
$f_{i+1}=\cdots =f_n=0$. Let $\Psi_{i+1}$ be the orientation on $N_{i+1}$ near $p$ such that the form
$$\Psi_{i+1}\wedge df_{i+1,1}\wedge df_{i+1,2}\wedge\cdots\wedge
df_{n,1}\wedge df_{n,2}$$ is positive for the given orientation on $M$. Note that the set $N_1$ is  discrete when
$\dim M=2n$ and $\Psi_1 =\pm 1$ in this case.

\noindent {\bf Theorem 3.3.} \ {\it The Chern classes $c_k(E_{{\bf C}})$ are Poincar\'{e} dual to the cycle
$D_{n-k+1}$ respectively, \ $k=1,\cdots,n$. Then $c_k(E_{{\bf C}})=0,$ if $D_{n-k+1} =D_{n-k}$ or $D_{n-k+1}$ is a
boundary. }

 \noindent{\it Proof}. \ The theorem has been proved in [6] by using the
Grassmann manifolds, see also [1], p.134. In the following we give a direct proof.

For any $k\geq 1$, we have $c_k(F_k)= \pi_{n-k}^*c_k(E_{{\bf C}})$. Then there is a $2k-1$ form $\eta_k$ on
$V(E_{{\bf C}},n-k+1)$ such that
$$\alpha_{n-k+1}^*c_k(\Omega_{F_k})=-d\eta_k.$$
 Restricting $\eta_k$ to
each fibre of $\alpha_{n-k+1}: \ V(E_{{\bf C}},n-k+1)\to V(E_{\bf C}, n-k)$ is the volume form. By Gram-Schmidt
process, from the sections $s_1,\cdots,s_{n-k+1}$, we have a section $\widetilde \sigma_{n-k+1}:
M-D_{n-k+1}(\sigma)\to V(E_{{\bf C}},n-k+1)$, where $\widetilde\sigma_{n-k+1}=\{\widetilde s_1,\cdots,\widetilde
s_{n-k+1}\}$ is a Hermite frame field on $M-D_{n-k+1}(\sigma)$. Let $U_i(\varepsilon)$ be a
$\varepsilon$-neighborhoods of $D_i$ respectively, $i=1,\cdots,n$, $U_i(\varepsilon)\subset U_j(\varepsilon)$ if
$i<j$. Then for any closed $m-2k$ form $\xi$, we have $$\int_M c_k(\Omega_{E_{{\bf C}}})\wedge \xi
=\lim\limits_{\varepsilon \to 0} \int_{\partial U_{n-k+1}(\varepsilon)}
\widetilde\sigma_{n-k+1}^*\eta_k\wedge\xi.$$

The space ${\partial U_{i}(\varepsilon)}- {\partial U_{i-1}(\varepsilon)}$ can be viewed as a fibre bundle over
 $N_i$ with  fibre  $S^{2n-2i+1}$.
The normal bundle $v(N_{i+1})$  of $N_{i+1}$ in $M$ is oriented and the orientation is given by that of $M$ and
$N_{i+1}$. Let $x=(x_1,\cdots,x_{m-2k},x_{m-2k+1},\cdots,x_m)$ be oriented coordinates in a neighborhood of $p\in
N_{n-k+1}$ in $M$ such that, restricting on $N_{n-k+1}$,  $(x_1,\cdots,x_{m-2k})$ be oriented coordinates on
$N_{n-k+1}$. By
\begin{eqnarray*} & & \Psi_{n-k+1}\wedge df_{n-k+1,1}\wedge
df_{n-k+1,2}\wedge\cdots\wedge df_{n,1}\wedge df_{n,2} \\
& &  =\frac {\partial (f_{n-k+1,1},f_{n-k+1,2},\cdots,f_{n,1},f_{n,2})}{\partial
(x_{m-2k+1},\cdots,x_m)}\Psi_{n-k+1}\wedge dx_{m-2k+1}\wedge\cdots\wedge dx_m, \end{eqnarray*} we have $\frac
{\partial (f_{n-k+1,1},f_{n-k+1,2},\cdots,f_{n,1},f_{n,2})}{\partial (x_{m-2k+1},\cdots,x_m)}> 0$. As noted above,
integrating along the fibres of the map
 $\partial U_{n-k+1}(\varepsilon)-\partial U_{n-k}(\varepsilon) \to N_{n-k+1}$  yields,
$$\lim\limits_{\varepsilon \to 0} \int_{\partial U_{n-k+1}
(\varepsilon)-\partial U_{n-k}(\varepsilon)} \widetilde\sigma_{n-k+1}^*\eta_k\wedge\xi =\int_{N_{n-k+1}}\xi.$$ As
$2n-2i+1 > 2k-1$ when $i\leq n-k$, we have
$$\lim\limits_{\varepsilon \to 0} \int_{\partial U_{n-k}
(\varepsilon)} \widetilde\sigma_{n-k+1}^*\eta_k\wedge\xi =0.$$ Hence
$$\int_M c_k(\Omega_{E_{{\bf C}}})\wedge \xi =\int_{D_{n-k+1}}\xi.$$

This completes the proof of the theorem. \ \  \ $\Box$

\noindent {\bf Corollary 3.4.} \ {\it If $N_{n-k+1}$ is a closed submanifold of $M$, that is, $\overline
{N_{n-k+1}}= N_{n-k+1}$. Then $N_{n-k+1}$ is the Poincar\'{e} dual of the Chern class $c_k(E_{{\bf C}})$.
Furthermore,  there is an oriented real vector bundle $F_{\bf R}$ over $N_{n-k+1}$ if $\dim_{\bf R} M = 2n$, such
that
$$\int_M c_k(\Omega_{E_{{\bf C}}})\wedge c_{n-k}(\Omega_{E_{{\bf C}}})
 =\int_{N_{n-k+1}} e(\Omega_{F_{\bf R}}),$$
 where $ \Omega_{F_{\bf R}}$  is the curvature of
$F_{\bf R}$ with respect to some connection on $F_{\bf R}$. }

 \noindent{\it Proof}. \ By Theorem 3.3, we have
$$\int_M c_k(\Omega_{E_{{\bf C}}})\wedge c_{n-k}(\Omega_{E_{{\bf C}}})
 =\int_{N_{n-k+1}} p^*c_{n-k}(\Omega_{E_{{\bf C}}}),$$
where $p: \ N_{n-k+1} \to M$  is the inclusion, $\dim N_{n-k+1} =2n-2k$. By dimensional reason, the pull-back
bundle $p^*E_{\bf C}$ can be decomposed by $p^*E_{{\bf C}} \cong F_{{\bf C}}\oplus {\cal E}_{k}$ on $N_{n-k+1}$,
where ${\cal E}_{k}$ is a trivial bundle. Hence
$$p^*c_{n-k}(E_{{\bf C}})=c_{n-k}(F_{{\bf C}}) = e(F_{\bf R}). \ \ \ \ \Box$$

\noindent {\bf Remark}. \ On the other hand, by the assumption of Corollary 3.4,  there is a nature decomposition
$p^*E_{{\bf C}} = \widetilde F_{{\bf C}}\oplus {\cal E}_{n-k}$, where ${\cal E}_{n-k}$ is a trivial bundle
generated by the sections $s_1,\cdots, s_{n-k}$ restricted on $N_{n-k+1}$. If $k<\frac n2$, we have
$$p^*c_{n-k}(E_{{\bf C}})= c_{n-k}(\widetilde F_{{\bf C}})=0.$$
Thus, if $\int_M c_k(\Omega_{E_{{\bf C}}})\wedge c_{n-k}(\Omega_{E_{{\bf C}}})\not= 0$, then there does not exist
generic sections  on the vector bundle $E_{{\bf C}}$ such that $N_{n-k+1}$ is a closed submanifold of $M$.

Notice that when $\dim M =2n$, $\int_Mc_n(E_{{\bf C}})$ is the intersection number of $s_1(M)$ with $M$ in
$TE_{{\bf C}}$, where the orientation on the fibres of $E_{{\bf C}}$ are determined by the complex structure.

\noindent {\bf  Theorem 3.5.} \ {\it Assuming $\dim M > 2n$, let $i: \ S \to M$ be an imbedding which intersect
transversally with $N_1$, where $S$ be a $2n$ dimensional oriented submanifold. Then $\int_S c_n(i^* E_{{\bf C}})$
is the intersection number of $S$ with $N_1$. }

\noindent{\it Proof}. \ It is easy to see that the section $s_1$ of the vector bundle $E_{{\bf C}}$ pulls back to
a section $i^*s_1$ of the bundle $i^*E_{{\bf C}}\to S$. The zeros of the section $i^*s_1$ correspond exactly to
the points of intersection of $S$ with $N_1$. If $p$ is a point in $S\cap N_1$, we have
$$T_pM=T_pN_1\oplus v_p(N_1) = T_pN_1\oplus T_pS,$$ where
$v_p(N_1)$ is the normal space of $N_1$ at $p$. Furthermore, the tangent map of $s_1: \ S \to E_{{\bf C}}$ at $p$
is an isomorphism of $T_pS$ to the fibre of $E_{{\bf C}}$ at $p$. The tangent map $s_{1*p}: \ T_pS \to E_{{\bf
C}}|_p$ preserving the orientation if and only if the orientation of $T_pS\otimes T_pN_1$ defined by that of
$T_pS$ and $T_pN_1$ is the same as that of $T_pM$, see the proof of Theorem 3.3. These completes the proof of the
theorem. \ \ \ $\Box$

As $p_k(E)=(-1)^kc_{2k}(E\otimes {\bf C})$, the Pontrjagin classes are Poincar\'{e} duals to some cycles on the
base manifold by Theorem 3.3. For example, let $TM$ be a tangent bundle of a $4k$-Riemannian manifold. Then there
is a decomposition: $TM\otimes {\bf C} = F_{2k}\oplus {\cal E}_{2k}$ and
$$p_k(TM) =(-1)^kc_{2k}(TM\otimes {\bf C})= (-1)^kc_{2k}(F_{2k}) =
(-1)^k e(F_{2k{\bf R}}).$$ Hence the Poincar\'{e} dual of the Pontrjagin class $p_k(TM)$ is  that of the Euler
class $(-1)^ke(F_{2k{\bf R}})$. In what follows we give some further study.

Let  $\sigma =\{s_1,\cdots,s_{n}\}$ be a set of sections of a real vector bundle $ E$ with the degeneracy set
$D_i=D_i(\sigma)$ defined as in the complex case. We call $\sigma$ generic if the sections $\sigma$ satisfies the
similar conditions. Denote $N_{i+1}= D_{i+1} -D_i, \  \dim N_{i+1}=m-n+i$, where $m=\dim M$.

\noindent {\bf Theorem 3.6.} \ {\it Let $\pi: \ E\to M$ be an oriented Riemannian vector bundle and $\sigma$ be a
set of generic sections defined as above.  If $N_{n-2k+1}$ is a closed submanifold of $M$, the Poincar\'{e} dual
of Pontrjagin class $p_k(E)$ can be represented by the Poincar\'{e} dual of Euler class of the normal bundle
$v(N_{n-2k+1})$ of $N_{n-2k+1}$ in $M$. Then $p_k(E)=0$ if the bundle $v(N_{n-2k+1})$ has a nowhere vanish section
or the cycle $D_{n-2k+1}$ is a boundary. }

\noindent{\it Proof}. \ For each $k$, we have $E|_{M-D_{n-2k}} = \hat F_{2k}\oplus {\cal E}_{n-2k}$, where ${\cal
E}_{n-2k}$ is generated by the sections $s_1,\cdots,s_{2k}$. The  bundle $\hat F_{2k}$ is oriented and
$s_{n-2k+1}$ is a transversal section of this bundle. In Proposition 12.7 of [1], Bott and Tu  proved that the
vector bundle $\hat F_{2k}\to N_{n-2k+1}$ is isomorphic to the normal bundle of $N_{n-2k+1}$ in $M$. These also
shows  the submanifold $N_{n-2k+1}$ is oriented.

By Proposition 2.4, we have $\pi_{n-2k}^*p_k(E) =[e(\widetilde F_{2k})]^2.$ There is a $2k-1$ form $\eta_k$ on
$V(E,n-2k+1)$ such that $\alpha_{n-2k+1}^*e(\Omega_{2k})=-d\eta_k,$ where $\Omega_{2k}$ is the curvature matrix of
the naturally defined connection on the bundle $\widetilde F_{2k}$. Form $s_1,\cdots,s_{n-2k+1}$, we have a
section  $\widetilde \sigma$ of $\pi_{n-2k+1}: \ V(E,n-2k+1)\to M$ on $M-D_{n-2k+1}$. It is easy to see that, on
$M-D_{n-2k+1}$,
$$p_k(E)  =\widetilde \sigma^*\pi_{n-2k+1}^*p_k(E) = -d\eta_k\wedge \widetilde
\sigma^*\alpha_{n-2k+1}^*e(\widetilde F_{2k}).$$ On the submanifold $N_{n-2k+1}$, we have
$$\widetilde
\sigma^*\alpha_{n-2k+1}^*e(\widetilde F_{2k})= (\alpha_{n-2k+1}\circ \widetilde \sigma)^*e(\widetilde
F_{2k})=e(\hat F_{2k})=e(v(N_{n-2k+1})).$$

The rest of the proof is similar to that of Theorem 3.3.
 Let $U_{n-2k+1}(\varepsilon)$ be a
$\varepsilon$-neighborhood of $D_{n-2k+1}$ in $M$. For any $m-4k$ form $\xi$ on $M$, we have
\begin{eqnarray*}
& & \int_M p_k(\Omega_E)\wedge \xi \\
  & & = -\lim\limits_{\varepsilon \to 0}
  \int_{M-U_{n-2k+1}(\varepsilon)} d\eta_k\wedge
  \widetilde\sigma^*\alpha_{n-2k+1}^*e(\Omega_{2k})\wedge\xi \\
& & =\int_{N_{n-2k+1}}\widetilde\sigma^*\alpha_{n-2k+1}^*e(\Omega_{2k})\wedge
 \xi,
 \end{eqnarray*}
where $\widetilde\sigma^*\alpha_{n-2k+1}^*e(\Omega_{2k})$ is the Euler form of the vector bundle $\hat F_{2k}\to
N_{n-2k+1}$. If the submanifold $N_{n-2k+1}$ is  closed,  then $\partial N_{n-2k+1}=\emptyset$ and
$$\int_{N_{n-2k+1}}\widetilde\sigma^*\alpha_{n-2k+1}^*e(\Omega_{2k})\wedge
 \xi =\int_{N_{n-2k+1}}e(\widetilde \Omega_{2k})\wedge
 \xi,$$
 where $\widetilde \Omega_{2k}$ is the curvature matrix of the
 normal bundle $v(N_{n-2k+1})$.
 The
theorem is proved. \ \ \ $\Box$

 In Proposition 12.8 of [1], Bott and Tu proved this kind
theorem for the Euler class.

\noindent {\bf Lemma 3.7.} \  {\it If $\dim M  = $ rank $E = 4k$, we can choose a set of generic sections such
that $\overline {N_{i}}\cap D_{i-1}=\emptyset$ for each $i\leq 2k+1$. Then $N_{i}$ are closed submanifolds of $M$
for $i\leq 2k+1$. }

\noindent {\it Proof}. \ We prove the lemma by induction. Assuming
 we have chosen generic sections $s_1,\cdots,s_{2k}$ such that
the lemma holds for each $i\leq 2k$. Then the set $N_{i}$ are all closed submanifolds of $M$ and $N_{i}\cap
N_{j}=\emptyset$ for any $i\not = j\leq 2k$. On $M-\bigcup N_{i}, E$ can be decomposed by $E|_{M- \bigcup N_{i}} =
\hat F_{2k} \oplus {\cal E}_{2k}$, where ${\cal E}_{2k}$ is  generated by $s_1,\cdots,s_{2k}$.  Let $U_i$ be a
neighborhood of $N_{i}$ such that $N_{i}$ is a deformation retract of $U_{i}$. Since $\dim N_{i}=i-1 < 2k$, the
bundle $\hat F_{2k}\to U_i-N_i$ has a nowhere vanish section. Then we can construct a section $\tilde s$ on
$E|_{U_i}$ such that $\tilde s$ is nowhere zero on $U_i-N_{i}$ and $\tilde s|_{N_{i}}=0$.
 By the partition of unit we have a section $\tilde s$ on $E$ such
that $s_1,\cdots,s_{2k},\tilde s$ are linearly independent on each $U_i - N_{i}$. With a perturbation of $\tilde
s$ on $M-\bigcup U_i$, we can get the desired section $s_{2k+1}$, cf. [1], p.123. \ \ \ $\Box$

\noindent {\bf Corollary 3.8.} \ {\it Let $M$ be an oriented manifold of dimension $4k$ and $N=N_{2k+1}$ be a
closed submanifold defined as in Lemma 3.7 for an oriented vector bundle $E$ of rank $4k$. Then we have
$$ \int_M p_k(E)= \chi (v(N)),$$
where $\chi(v(N))$ is the Euler characteristic of the normal bundle $v(N)$. }

Let $TM$ be the tangent bundle of a 4-dimensional oriented manifold. Then $\int_M \frac 13p_1(M)=\frac 13\chi
(v(N))$ is the signature of the manifold $M$. As we know, the signature of 4-manifold $M$ is a multiple of $8$.
Furthermore, if $M$ is spin, sig$(M)$ is a multiple of $16$, see for example [8], p.280. Then $\chi (v(N))$ is a
multiple of $24$ or $48$ respectively.

 \vskip 1cm \centerline{\bf  References} \vskip 0.3cm

{\small \noindent [1] \  R. Bott  and L. Tu,  Differential forms in algebraic topology, \ Springer GTM 82, 1982.

\noindent [2] \ S. S. Chern, A simple intrinsic proof of the Gauss-Bonnet formula for the closed Riemannian
manifolds, Ann. of Math., {\bf 45}(1944), 747-752.

\noindent [3] \ S. S. Chern, On the curvature integral in a Riemannian manifold, Ann. of Math., {\bf 46}(1945),
674-648.

\noindent [4] \ S. S. Chern, Characteristic classes of Riemannian manifolds, Ann. of Math., {\bf 47}(1946),
85-121.

\noindent [5] \ S. S. Chern, On curvature and characteristic classes of a Riemannian manifold, Abh. Math. Sem.
Univ. Hamburg, {\bf 20}(1955), 117-162.

\noindent [6] \ P. Griffiths  and J. Harris, Principles of algebraic geometry, Wiley-Interscience, New York, 1978.

\noindent [7]  \ S. Kobayashi and K. Nomizu, \ Foundations of differential geometry, vol. 2,  Interscience
Publishers, New York, 1969.

\noindent    [8] \ H. B. Lawson  and H. Michelsohn, Spin geometry, Princeton University Press, Princeton, New
Jersey, 1989.

\noindent [9] \ V. Mathai  and D. Quillen, Superconnections, Thom classes and equivariant differential forms,
Topology, {\bf 25}(1986), 85-110.

\noindent [10] \ J. W. Milnor and J. D. Stasheff, Characteristic classes, Ann. of Math. Studies, no. 76, Princeton
University Press, 1974.  

 \vskip 0.5cm

 \noindent E-mail: \ jwzhou@suda.edu.cn

\end{document}